**Forward period analysis and the long term simulation of a periodic Hamiltonian system**

Wang Pengfei[1,2]

[1]Center for Monsoon System Research, Institute of Atmospheric Physics, Chinese Academy of Sciences, Beijing 100190

[2]State Key Laboratory of Numerical Modeling for Atmospheric Sciences and Geophysical Fluid Dynamics, Institute of Atmospheric Physics, Chinese Academy of Sciences, Beijing 100029

**Abstract**

The period of a Morse oscillator and mathematical pendulum system are obtained, accurate to 100 significant digits, by forward period analysis (FPA). From these results, the long-term [0, $10^{60}$] (time unit) solutions, which overlap from the Planck time to the age of the universe, are computed reliably and quickly with a parallel multiple-precision Taylor series (PMT) scheme. The application of FPA to periodic systems can reduce the computation loops of long-term reliable simulation from $O\left(t^{1+1/M}\right)$ to $O\left(\ln t + T/h_0\right)$ where $T$ is the period, $M$ the order and $h_0$ a constant step-size. This scheme provides a way to generate reference solutions to test other schemes' long-term simulations.

**1 Introduction**

The Hamiltonian system is important for many areas of physical research, and the Hamiltonian constant, $H$, is a characteristic quantity of the system. Earlier numerical methods, such as the Euler, Runge–Kutta and linear multistep methods, as well as some low-order Taylor methods, are non-structure-preserving methods because they cause the structure of the system to change so that *H* varies with time, *t*. When applied to simulate a Hamiltonian system, they can induce unstable computation or lead to incorrect solutions. The symplectic method, developed by a number of researchers [1-3], is a common structure-preserving method. It conserves the area or volume of the system during computation. The square conservation also preserves structure, through conserving the length of the simulated system [4]. These structure-preserving methods allow the Hamiltonian constant to remain almost constant, or only change periodically, during the entire integration time range $[0, t]$. The key advantages of structure-preserving methods are that they provide a true long-term trajectory of the simulated system and stabilize the computation process.

However, these structure-preserving methods still have minor problems. First, when dealing with general nonlinear Hamiltonian systems, some of the symplectic methods are based on the generating function. These implicit methods are applied to solve nonlinear algebraic equations at each step, and thus efficiency becomes a problem. Some symplectic methods are based on the Runge–Kutta method, the order of which is generally under 10 (and rarely above 15), to avoid a complicated procedure. Some other high-order explicit symplectic methods are used to study separable Hamiltonian systems [2, 5], but these explicit methods usually still have an order of less than 10, and are limited to separable Hamiltonian systems.

Second, although an advantage of symplectic methods is that the step-size can be chosen with a larger value than in a classical approach – saving simulation time – such an increase in step-size causes large errors in the primary variables. This happens despite no modification of the trajectory structure. Gladman [6] has said that "…the conservation of the integrals is not a problem for the SIAs (symplectic integration algorithms) but the phase errors can still be uncomfortable after a large number of orbits. If one wanted to integrate the Earth for the lifetime of the solar system it is doubtful that these two SIAs could perform the $\sim 10^9$ orbit integration reliably. This is not necessarily too disheartening, since no other integration scheme known to the authors could perform the integration either."

The phase trajectory of a Hamiltonian system is one of the most basic requirements, but considering the period of the Hamiltonian system, the symplectic method provides no special insight, and only gives approximate numerical periods in which precision is proportional to the order of method and step-size. The long-term integration for a dynamical system is a challenging but necessary task in many subjects. Orders of magnitude of time periods in physics range from the Planck time ($5.39 \times 10^{-44}$ s) to the age of universe ($\sim 10^{17}$ s). Thus, a meaningful non-dimensional time for a dynamical system is within $10^{60}$ orders of magnitude. Simulation of the position of dynamical variables (not the trajectory structure) in the ultra long-term (i.e. $t = 10^{60}$) is still a time consuming task, even for a periodic dynamical system that applies the symplectic method or other existing methods. In this study, the author presents a procedure to solve such issues.

**2 Direct simulation of a Hamiltonian system and preservation of the total energy by the parallel multiple-precision Taylor method**

The parallel multiple-precision Taylor (PMT) method [7-12] was originally designed to solve nonlinear chaotic systems. It can provide ultra-high reliable solutions for longer times than other methods. The order of the PMT method can be very high compared to traditional approaches. Here, the application of the PMT method to a nonlinear Hamiltonian system is demonstrated. The orbits of the system for two atoms with Morse potential energy [13] are

$$\begin{cases} \dfrac{dx}{dt} = p \\ \dfrac{dp}{dt} = e^{-2x} - e^{-x} \end{cases} \text{ and } H = \dfrac{p^2}{2} + \dfrac{1}{2}\left(e^{-2x} - 2e^{-x}\right). \tag{1}$$

The initial values are $x_0 = 0$ and $p_0 = \sqrt{1-0.02}$ where $x$ represents displacement and $p$ momentum of the particles. Using the substitution $y = e^{-x}$ yields

$$\begin{cases} \dfrac{dy}{dt} = -py \\ \dfrac{dp}{dt} = y^2 - y \end{cases}, \tag{2}$$

where $y_0 = 1$ and $p_0 = \sqrt{1-0.02}$, and $H = \dfrac{p^2}{2} + \dfrac{1}{2}\left(y^2 - 2y\right)$ is still a constant.

The Taylor series expansions relevant to solving (2) are $\begin{cases} y_{n+1} = y_n + \sum\limits_{k=1}^{M} \alpha_k h^k \\ p_{n+1} = p_n + \sum\limits_{k=1}^{M} \beta_k h^k \end{cases}$, where the coefficients are given by $\alpha_k = \dfrac{1}{k!}\dfrac{d^k y(t_n)}{dt^k}$ and $\beta_k = \dfrac{1}{k!}\dfrac{d^k p(t_n)}{dt^k}$; and $h$ is the step-size. The coefficients are determined from the initial conditions $\begin{cases} \alpha_0 = y_n \\ \beta_0 = p_n \end{cases}$, and the relations

$$\begin{cases} \alpha_{k+1} = \dfrac{-1}{k+1}\left(\sum\limits_{i=0}^{k} \alpha_{k-i}\beta_i\right) \\ \beta_{k+1} = \dfrac{1}{k+1}\left(\sum\limits_{i=0}^{k} \alpha_{k-i}\alpha_i - \alpha_k\right) \end{cases}.$$ More detail regarding the parallel scheme can be found in Wang et al. [7]. The solution of (1) can be expressed as $x_{n+1} = -\ln y_{n+1}$, after obtaining the values of $y_{n+1}$ and $p_{n+1}$. If the order, $M$, is larger than 100, the parallel scheme greatly

decreases the computation time; if $M$ is smaller than 30, generally one CPU is sufficient to perform the calculation on a reasonable timescale.

The symplectic method used to solve (1) is a 2$^{nd}$ order explicit method (SE2), introduced by Qin et al. [14]:

$$u_1 = p_k - hc_1 f(x_k), v_1 = x_k + hd_1 g(u_1),$$

$$p_{k+1} = u_1 - hc_2 f(v_1), x_{k+1} = v_1 + hd_2 g(p_{k+1}),$$

where $c_1 = 0$, $c_2 = 1$, $d_1 = d_2 = \frac{1}{2}$, $f(x) = -(e^{-2x} - e^{-x})$, and $g(p) = p$.

The second and more complex example system is a non-separable Hamiltonian system, defined by:

$$\begin{cases} \dfrac{dp}{dt} = p \sin q \\ \dfrac{dq}{dt} = p + \cos q \end{cases}. \tag{3}$$

The Hamiltonian is $H = \dfrac{1}{2} p^2 + p \cos q$, and the initial values are $q_0 = 0, p_0 = 1$:

$$\begin{cases} p_{n+1} = p_n + \sum_{k=1}^{M} \alpha_k h^k \\ q_{n+1} = q_n + \sum_{k=1}^{M} \beta_k h^k \end{cases},$$

where $b(t) = \sin q$, $g(t) = \cos q$, $c(t) = p \sin q$, $d(t) = p + \cos q$ and the $k^{th}$ Taylor coefficients are $b_k$, $g_k$, $c_k$, and $d_k$. Therefore, other coefficients are: $\alpha_{m+1} = \dfrac{1}{m+1} c_m$, $\beta_{m+1} = \dfrac{1}{m+1} d_m$, $b_m = \dfrac{1}{m} \sum_{i=1}^{m} i g_{m-i} \beta_i$, $g_m = \dfrac{-1}{m} \sum_{i=1}^{m} i b_{m-i} \beta_i$, $c_m = \sum_{i=0}^{m} \alpha_{m-i} b_i$, $d_m = \alpha_m + g_m$, and the initial coefficients are $\alpha_0 = p_0$, $\beta_0 = q_0$, $b_0 = \sin q_0$, $g_0 = \cos q_0$, $c_0 = p_0 \sin q_0$, and $d_0 = p_0 + \cos q_0$.

The conservation of the Hamiltonian implies $\frac{dH}{dt}=0$. In the present numerical simulation, if $|\Delta H|\leq 10^{-16}$ (i.e. the smallest relative error in double-precision floating point arithmetic;) it is regarded as unchanged. When $H$ is a non-zero constant, $\left|\frac{\Delta H}{H}\right|\leq 10^{-16}$ is a criterion for a large Hamiltonian.

For instance, there is a Hamiltonian, $H=-0.01$, for Equation (2) when the numerical solutions of $p$ and $x$ are $p_N$ and $x_N$. The error of the Hamiltonian is

$$\Delta H = H_N - H = \left\{\frac{p_N^2}{2}+\frac{1}{2}(y_N^2-2y_N)\right\}-\left\{\frac{p^2}{2}+\frac{1}{2}(y^2-2y)\right\}.$$

Since $|\Delta H|\leq\left|\frac{p_N^2}{2}-\frac{p^2}{2}\right|+\frac{1}{2}|y_N^2-y^2|+|y_N-y|$, $\left|\frac{p_N^2}{2}-\frac{p^2}{2}\right|+\frac{1}{2}|y_N^2-y^2|+|y_N-y|<10^{-16}$ guarantees $|\Delta H|<10^{-16}$. The numerical error at time $t$ indicates that $|p_N-p|\leq C_1 h^{M+1}$, $|y_N-y|\leq C_2 h^{M+1}$, $\left|\frac{p_N^2}{2}-\frac{p^2}{2}\right|\approx|p(p_N-p)|\leq|p|C_1 h^{M+1}$ where $C_1$ and $C_2$ are constants, and $\left|\frac{y_N^2}{2}-\frac{y^2}{2}\right|\approx|y(y_N-y)|\leq|y|C_2 h^{M+1}$. Since $|p|$ and $|y|$ are bounded variables, $\left|\frac{p_N^2}{2}-\frac{p^2}{2}\right|+\frac{1}{2}|y_N^2-y^2|+|y_N-y|<Ch^{M+1}$, where $C$ is a constant that satisfies $|p|C_1+|y|C_2+C_2\leq C$.

With a step-size of $h=0.01$, a high enough order $M$ can always be chosen to guarantee $\left|\frac{p_N^2}{2}-\frac{p^2}{2}\right|+\frac{1}{2}|y_N^2-y^2|+|y_N-y|<Ch^{M+1}<10^{-16}$. In practice, the order of $M$ can be easily determined by several numerical runs without knowing the value of $C$. By using this high-order method, the structure-preserving solution of the original equation is obtained by numerical means. In fact, because increasing $M$ is very easy to do with the Taylor series method, we can make $|\Delta H|$ even smaller for Equation (2).

In the direct simulation of Eq. (2) with $t = 10^7$, a 20-order PMT scheme was used to achieve the simulation results. From Fig. 1c, the PMT method is shown to predict the correct trajectory structure (*x-p* plane), and the cycle of *x* is also correct (Fig. 1a). During the entire computation time range, the Hamiltonian *H* approaches a constant (Fig. 1d), while Fig. 1e indicates that *H* varies periodically and has larger errors when simulated with SE2. The more important issue is that the error in *x* increases as the simulation time increases. Thus, the position of *x* is not reliable at times longer than $10^5$ time units (Fig. 1b, see also Table 4). Meanwhile, applying the PMT method to solve the non-separable Eq. (3), the variable (Fig. 2b) and the Hamiltonian *H* are simulated well, and *H* remains a constant throughout the simulation time range.

The PMT scheme is a self-verifying scheme, as discussed in [7]. This verification scheme is standard operation for PMT and CNS [12, 15, 16] numerical experiments.

## 3 Forward period analysis to obtain the period of a Hamiltonian system

Establishing the phase trajectories of Hamiltonian systems is a basic requirement, and the symplectic method achieves this as well as PMT. However, the symplectic method only gives approximate numerical periods, with a precision proportional to the order of the method and step-size.

The key of numerical method to identify a period of a dynamical system is to find out the solutions return to the initial valuess (if the system defined by *n*-th 1[rst] order differential equation, the *n* variables must return to their initial values simultaneously). The integration time between the start point and the repeat point thus approximate to the period. Generally the period obtained by numerical method in this way varies. The error between a variable (such as $x$) and its corresponding repeat point (regard as $E_x$, for example $E_x \leq 10^{-30}$) indicates the uncertainty ($\Delta t$) in estimating the period. Because $E_x$ is small and $\frac{dx(t)}{dt} \approx \frac{E_x}{\Delta t}$, $\Delta t \cdot p \approx E_x$. The standard division of $p$ was obtained through numerical experiments, $\sigma(p) = \sqrt{\overline{p^2}}$, as an averaged value of $p$ such that $\Delta t \sim E_x / \sigma(p)$. This formula suggests that the error bounds of a typical period have a magnitude of $10^{-30}$.

The forward period analysis (FPA) method is proposed to obtain the period for Eq. (2). The first stage of FPA is a pre-computation to find a suitable residual interval. The computation starts with $y = 1$, $p = \sqrt{1-0.02}$ and $\dot{y} = -py = -\sqrt{1-0.02}$, and a time-step size of $h_0 = 0.01$. At each step the computed values of $y_k$ and $y_{k-1}$ are checked to find out the approximate first period. The first repeat position satisfies $y \approx 1$ and $\dot{y} \approx -\sqrt{1-0.02}$, determined within the interval of $[T_l, T_h]$, where $T_l$ is the lower and $T_h$ the higher bound. This interval is the time it takes for $y$ to cross the base line $y = 1$ (from the $y > 1$ to the $y < 1$ direction, the steps that reach the lower bound are regarded as $k$). The first period is now between $T_l = kh_0$ and $T_h = (k+1)h_0$ time units, and thus the period $T \approx T_l$. The error of the period is about $(h_0 + \Delta t) \approx 10^{-2}$, and thus the precision of this forward period analysis method is mostly dependent on the last computation step-size $h_F$ (in this stage, $h_F = h_0$).

The second stage is the post-computation at the residual interval $[kh_0, (k+1)h_0]$ where $k = \text{Int}(t/T)$ is a positive integer number. This interval can be separated into subintervals by the dichotomy method. Denoting the whole interval as before, $[T_l, T_h]$, with $T_l = kh_0$ and $T_h = (k+1)h_0$, a new step size $h_F = h_F / 2$ is chosen to separate the interval into $[T_l, T_l + h_F]$ and $[T_l + h_F, T_h]$. If the value of $y$ at $T_l + h_F$ does not cross the base line then $T_l = T_l + h_F$; otherwise, $T_h = T_l + h_F$ at program and then repeat the operation in the new interval $[T_l, T_h]$. The dichotomy method maintains $h_F$ smaller than the magnitude of $E_T / \sigma(p)$, and thus total error of the period is dominated by $\Delta T$ ($\Delta T < T_h - T_l \approx 10^{-30}$). The computation cost for the dichotomy method in the last interval is about $30\log_2 10 \approx 100$ loops, while the pre-computation stage being about $T/0.01$ loops. The value of $T$ can be roughly estimated from Fig. 3a ($T$ is within a 45 time unit). The total loops are within 5000 loops for obtaining the period of Eq. (2).

The above procedures are also suitable for SE2 and other symplectic methods; if we choose the step-size ($h$) for SE2, the period which precise at $h^2$ magnitude is obtained. Applying the dichotomy method at the last interval for SE2, the error at $T_l$ and $T_h$ must first be confirmed to small enough, and this is not the superiority for SE2. Since decreasing $h$ greatly increases the computation time, if a more accurate period is required, for example $\Delta T = 10^{-30}$, $h \approx 10^{-15}$ must be set, and this requires $10^{15}$ loops to finish the computation.

In addition, applying self-verification by decreasing step-size requires the step-size to be about $h/100$ or less in SE2 to guarantee the reference solution is more accurate than the solution which step-size is $h$, and this requires 100 times more loops than the computation process. While self-verification of the PMT method only requires increasing the order $M$ to ~$M$+10, this does not increase the number of computation loops of the verification process, only time cost per loop. The increasing time cost in one integration loop is insignificant when $M$ is large (for example $M$>100). As a consequence, the PMT method is efficiently verified.

From the enlarged time axis in Fig. 3b, the first stage of FPA is to determine the residual interval of Eq. (2) as $[44.42, 44.43]$, i.e. $T_l = 44.42$ and $T_h = 44.43$. The FPA procedure in this interval, and the values of to 30 and 100 significant digits, are listed in Table 1. To the author's knowledge, this level of accuracy has never been reported.

**Table 1.** The period of a Morse system obtained by FPA with $M = 200$, $h = 0.01$, and the precision we use is 2000 bits.

| Significant digits | $T$ |
|---|---|
| 30 | 44.4288293815836624701588099006 |
| 100 | 44.4288293815836624701588099006069369861462168937569022308539560 6956434793099473910575326934764765237 |

**4 The application of FPA in long-term simulations**

Before demonstrating an ultra long-term simulation of Eq. (2), a simple periodic dynamical system is analyzed to determine the most important parameters in the long-term computation. The simple dynamical system is defined by $\begin{cases} dx/dt = p \\ dp/dt = -x \end{cases}$ and the initial values are $\begin{cases} x(0) = 0 \\ p(0) = 1 \end{cases}$. The issue is how to obtain 16 significant digits of $x(t)$ at $t = 10^{30}$. Since the analytical solution of this equation is $\begin{cases} x(t) = \sin(t) \\ p(t) = \cos(t) \end{cases}$, the result should be $x(10^{30}) = \sin(10^{30})$. As $\sin(t - 2\pi k) = \sin(t)$, the result is given by $\sin(10^{30}) = \sin(10^{30} - 2\pi k)$, and $k = \text{Int}\left[\dfrac{10^{30}}{2\pi}\right]$, i.e. the integer part of $\dfrac{10^{30}}{2\pi}$. The precision of $\sin(10^{30})$ is dependent on the precision of the approximation to $2\pi$. The 50 significant digits reference value of $\pi$ and the computed $k$ are listed in Table 2, and the double-precision (16 significant digits) results are also compared. From Table 2, note that the $k$ correspondence to the two different precisions of $\pi$ are different. The different $k$-values cause the residual of $t$, i.e. $T_r = 10^{30} - 2\pi k$, to be more uncertain. Therefore, precision to 16 significant digits for $\sin(10^{30})$ is not possible in a double-precision float platform.

The above example indicates that the reliable long-term computation of a periodic system is dependent on the precision of the period; $2\pi$ can be regarded as the period in this example. The relative error bound, $\varepsilon$, is estimated for the period to guarantee the computation error at $t$ is limited such that $\Delta x = 10^{-16}$. The true residual time can be written as

$$T_r = 10^{30} - 2\pi \text{Int}\left[\frac{10^{30}}{2\pi}\right],$$

and the numerical error induced residual time is

$$T_r' = 10^{30} - 2\pi(1+\varepsilon)\text{Int}\left[\frac{10^{30}}{2\pi(1+\varepsilon)}\right].$$

The first restriction of $\varepsilon$ guarantees that $k = \text{Int}\left[\frac{10^{30}}{2\pi}\right]$ and $k = \text{Int}\left[\frac{10^{30}}{2\pi(1+\varepsilon)}\right]$ are the same. Under this situation, we have $T_r = 10^{30} - 2\pi k$, $T_r' = 10^{30} - 2\pi(1+\varepsilon)k$ and the difference between $T_r$ and $T_r'$ is $\Delta t = 2\pi k \varepsilon$. The error bounds, $\varepsilon$, satisfy $2\pi k\varepsilon < 10^{-16}$ will guarantee $\Delta x < 10^{-16}$. Since $0 < 10^{30} - 2\pi k < 2\pi$, $2\pi k \approx 10^{30}$ and $\varepsilon < 10^{-16}/2\pi k = 10^{-46}$ is the relative error bound of $2\pi$. The 50 significant digits of $\pi$ satisfy this error bound. Thus, $T_r = 3.231831977487846$ and $\sin(10^{30}) = -0.090116901912138058$.

**Table 2.** The values of $\pi$ and $k$.

| Significant digits | $\pi$ | $k$ |
|---|---|---|
| 50 | 3.14159265358979323846264338327 95028841971693993751 | 159154943091895335768883763372 |
| 16 | 3.141592653589793 | 159154943091895335768883763373 |

The analysis of error bounds for a general periodic dynamical system is the same as the example system by alternating the period from $2\pi$ to $T$. Thus, $\varepsilon < E(t)/t$ is the period error bound, where $t$ is the simulation time and $E(t)$ the required relative error bound for the output. The fast computation of $\sin(10^{30})$ is a benefit from the pre-known of precise value of $\pi$. However, for a general periodic dynamical system such as Eq. (2) there is no such pre-knowledge for the period $T$, hence proposing FPA to obtain the precise $T$ first.

Long-term simulation by FPA is achieved by dividing the long-term ($t \gg T$) computation into two parts: one is the period detection of a cycle; the other is the simulation of the residual time, equal to $t - kT$. Because the computation error of the symplectic method generally propagates linearly with time [17, 18], such that an increase of $t$ times requires a $t^{1/M}$ times smaller step-size to control the error, where $M$ denotes the order of the symplectic method. The computation complexity for time $t$ in unit loops is $O(t^{1+1/M})$. However, applying the FPA procedure with PMT can help to reduce the computation from $O(t^{1+1/M})$ to $O(T/h_0)$ in the first stage, and to $O(\ln t)$ in the second stage, i.e. $O(\ln t + T/h_0)$.

**Table 3.** The loops used in the different schemes.

| Scheme | Direct integration | Applying FPA |
|---|---|---|
| Symplectic | $O(t^{1+1/M})$ | $O(\ln t + Tt^{1/M})$ |
| PMT | $O(t/h_0)$ | $O(\ln t + T/h_0)$ |

As illustrated in Table 3, the long-term computation of a dynamical system uses at least $O(t)$ loops without FPA, but with FPA this is cut to $O(\ln t + T/h_0)$ at most. This improvement greatly decreases the CPU time cost and makes many unsolvable long-term problems reliably solvable.

**Table 4.** The variable *p* obtained by FPA and direct simulation for a Morse system.

| t | Direct integration with SE2 | Direct integration with PMT | By FPA with PMT |
|---|---|---|---|
| 0 | 0.989949493661166 | 0.989949493661166 | 0.989949493661166 |
| 10 | 0.142033683767425 | 0.142049967327890 | 0.142049967327890 |
| $10^2$ | 0.120638240019144 | 0.120968440888269 | 0.120968440888269 |
| $10^3$ | −0.015582896828410 | −0.013519353495639 | −0.013519353495639 |
| $10^4$ | 0.275552330918520 | 0.406695207104251 | 0.406695207104251 |
| $10^5$ | 0.120371703265026 | −0.209442226126745 | −0.209442226126745 |
| $10^6$ | −0.017153824006555 | −0.575071021680786 | −0.575071021680786 |
| $10^7$ | 0.214900072283681 | 0.406850634713920 | 0.406850634713920 |
| $10^8$ | - | - | −0.208888390114776 |
| $10^9$ | - | - | −0.548460247715134 |
| $10^{10}$ | - | - | 0.632436322896067 |
| $10^{20}$ | - | - | −0.459751580833174 |
| $10^{30}$ | - | - | 0.203324673559885 |
| $10^{40}$ | - | - | 0.240998707662065 |
| $10^{50}$ | - | - | 0.544244466975686 |
| $10^{60}$ | - | - | −0.839008449972302 |

Given the period obtained by FPA (see Table 1), the variable values for the Morse system can be computed quickly. The results of selected times from (10 to $10^{60}$) are listed in Table 4, with the corresponding values obtained by direct integration. The results of FPA and PMT direct integration agree well. The results of the SE2 method have errors from the early integration stage, and the results beyond $10^5$ are incorrect. It took about one day to obtain the result at $t=10^7$ by direct integration, so obtaining a result at $10^{60}$ is a seemingly impossible task for direct integration, but with the help of FPA, reliable results can be obtained.

Another classical dynamical system is the mathematical pendulum. The Hamiltonian of a pendulum system is $H = \frac{1}{2}p^2 - \cos q$, and the dynamical equation is $\begin{cases} \dot{q} = p \\ \dot{p} = -\sin q \end{cases}$ with initial values $\begin{cases} q(0) = 0 \\ p(0) = 1 \end{cases}$. The period of this system approaches to $2\pi$, while the initial momentum $p \to 0$ [17]. Table 5 lists the period corresponding to the initial condition, $q = 0$, with different momenta, $p$. All periods are accurate to 50 significant digits. As illustrated in Table 5, the period approaches $2\pi$ when $p$ decreases from 1 to $10^{-30}$. This experiment again proves the correctness of FPA.

**Table 5.** The period of a mathematical pendulum system obtained by FPA with $M = 200$, $h = 0.01$, and the precision we use is 2000 bits.

| $P$ | $T$ (50 significant digits) |
|---|---|
| 1 | 6.7430014192503841714848146311963079580032035765643 |
| $10^{-1}$ | 6.2871178299331781141446745665180361610970124356918 |
| $10^{-2}$ | 6.2832245776399990205430348375448192284997546476674 |
| $10^{-3}$ | 6.2831856998787233989673928392330685974085610706684 |
| $10^{-4}$ | 6.2831853111065772994348591606129983671149348353934 |
| $10^{-5}$ | 6.2831853072188563850957114151234372289079555651010 |
| $10^{-10}$ | 6.2831853071795864769292137573759930099424226253101 |
| $10^{-20}$ | 6.2831853071795864769252867665590057683943780686583 |
| $10^{-30}$ | 6.2831853071795864769252867665590057683943387987502 |
| $2\pi$ | 6.2831853071795864769252867665590057683943387987502 |

The speed of applying FPA to obtain a period for these demonstration systems is very fast, and the computation finished within 1 minute on a Linux system with an Intel Xeon 2.5 Ghz CPU. The long-time scope solutions obtained by FPA also finished within 1 minute.

**5 Discussion and conclusion**

Using FPA, the periods of some classical Hamiltonian systems are successfully obtained, accurate to 100 significant digits. The $t \in [0,10^{60}]$ reliable solution computation method and results are demonstrated. To the best of the author's knowledge, this accuracy of time period for long-term solutions of Hamiltonian systems has not been reported before. The FPA method provides a powerful tool to gain time-effective ultra long-term reliable solutions of periodic systems.

The FPA procedure works well in conjunction with the traditional symplectic method and the PMT method. Generally, symplectic methods with different orders require different subroutines to conduct the computation. In contrast, it is relatively easy to change the order of the Taylor series method, which is an advantage in that it has the flexibility to carry out simulations with different orders of accuracy for one system. The PMT method can directly simulate the Morse system well for $t \in [0, 10^7]$, but for much longer times simulation it is hard without FPA. The demonstration systems here are simple. For more complex systems, higher order PMT approaches can be used. Indeed, details of an example application of a high-order PMT method to directly simulate the Henon–Heiles system can be found in Liao [15].

The Taylor series method has a good convergence property when the order is high enough [9]. This feature can enlarge the step-size $h$ to 0.01 for Eq. (2), and increase the simulation speed. The result obtained by the Taylor series method not only maintains $\Delta H \simeq 0$, but also reduces numerical errors. The PMT method is not a structure-preserving method, but it can preserve the structure well by numerical means. Consequently, it can be used as an alternative to symplectic methods for the computation of simple Hamiltonian systems. Moreover, PMT can be applied to some non-separable nonlinear Hamiltonian systems as well as separable ones, and even to non-Hamiltonian chaotic systems.

The essence of applying FPA to long-term computation is divided into two parts: one is the period detection of a unit cycle; the other part is the computation of residue time equal to $t - kT$. This procedure helps to reduce the computation time for the long-term reliable simulation from $O\left(t^{1+1/M}\right)$ to $O\left(\ln t + T/h_0\right)$. The FPA procedure benefits not only the PMT method, but also the traditional symplectic method. The main problem with the symplectic method is that if the order *M* is not large enough (for example *M*<10) it still requires many computation loops for $t = 10^{60}$ – about $O\left(60\ln 10 + T10^{60/M}\right)$ loops, this takes a long time. For a medium-term time period, such as the ~$10^9$ orbits of the Earth–solar system, the solution is required at a $t = 10^{17}$ s magnitude. In this case, the symplectic method should work well with FPA.

## Acknowledgements

This research was jointly supported by the National Natural Sciences Foundation of China (41375112) and the National Basic Research Program of China (2011CB309704).

## References

1. Feng, K., *On difference schemes and sympletic geometry.* Feng Kang (Ed.), Proceedings of the 1984 Beijing Symposium on Differential Geometry and Differential Equations, Science Press, Beijing (1985), 1984: p. 42-58.
2. Ruth, R.D., *A canonical integration technique.* IEEE Trans. Nucl. Sci, 1983. **30**(4): p. 2669-2671.
3. Sanz-Serna, J., *Runge-Kutta schemes for Hamiltonian systems.* BIT Numerical Mathematics, 1988. **28**(4): p. 877-883.
4. Bin, W., Z.Q. Cun, and J. ZhongZhen, *Square conservation systems and Hamiltonian systems.* Science in China, Ser. A, 1995. **38**(10): p. 1211-1219.
5. Yoshida, H., *Construction of higher order symplectic integrators.* Physics Letters A, 1990. **150**(5): p. 262-268.
6. Gladman, B., M. Duncan, and J. Candy, *Symplectic integrators for long-term integrations in celestial mechanics.* Celestial Mechanics and Dynamical Astronomy, 1991. **52**(3): p. 221-240.
7. Wang, P.F., J.P. Li, and Q. Li, *Computational uncertainty and the application of a high-performance multiple precision scheme to obtaining the correct reference solution of Lorenz equations.* Numerical Algorithms, 2012. **59**(1): p. 147-159.
8. Liao, S.J., *On the reliability of computed chaotic solutions of non-linear differential equations.* Tellus A, 2009. **61**(4): p. 550-564.
9. Barrio, R., *Performance of the Taylor series method for ODEs/DAEs.* Applied Mathematics and Computation, 2005. **163**(2): p. 525-545.
10. Moore, R.E., *Interval analysis.* Vol. 2. 1966: Prentice-Hall Englewood Cliffs.

11. Wang, P., Y. Liu, and J. Li, *Clean numerical simulation for some chaotic systems using the parallel multiple-precision Taylor scheme.* Chinese Science Bulletin, 2014. **59**(33): p. 4465-4472.
12. Liao, S. and P. Wang, *On the mathematically reliable long-term simulation of chaotic solutions of Lorenz equation in the interval [0, 10000].* Sci China: Physics, Mechanics and Astronomy., 2014. **57**(2): p. 330-335.
13. Feng, K. and M. Qin, *Symplectic geometric algorithms for Hamiltonian systems.* 2010: Springer.
14. Qin, M., D. Wang, and M. Zhang, *Explicit symplectic difference schemes for separable Hamiltonian systems.* J. Comput. Math, 1991. **9**(3): p. 211-221.
15. Liao, S., *On the numerical simulation of propagation of micro-level inherent uncertainty for chaotic dynamic systems.* Chaos, Solitons & Fractals, 2013. **47**(1): p. 1-12.
16. Liao, S., *Physical limit of prediction for chaotic motion of three-body problem.* Communications in Nonlinear Science and Numerical Simulation, 2014. **19**(3): p. 601-616.
17. Hairer, E., C. Lubich, and G. Wanner, *Geometric numerical integration (Second Edition),* 2005, Springer. p. 636pp.
18. Quispel, G. and C. Dyt, *Volume-preserving integrators have linear error growth.* Physics Letters A, 1998. **242**(1): p. 25-30.

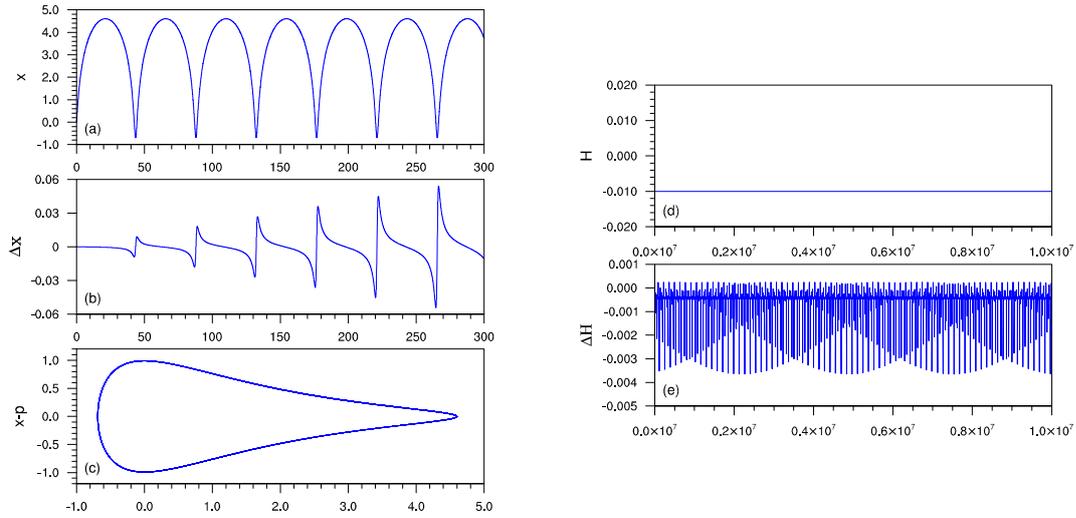

**Figure 1.** Illustrations of the direct simulation of Eq. (2): (a) variable *x* by PMT; (b) error of *x* computed by the SE2 method; (c) structure of the *x-p* plane by PMT; (d) the Hamiltonian *H* by PMT to $T = 10^7$; (e) error of the Hamiltonian *H* by the SE2 method. Panels (a–c) have step-size $h = 0.01$; (d–e) have step-size $h = 0.1$.

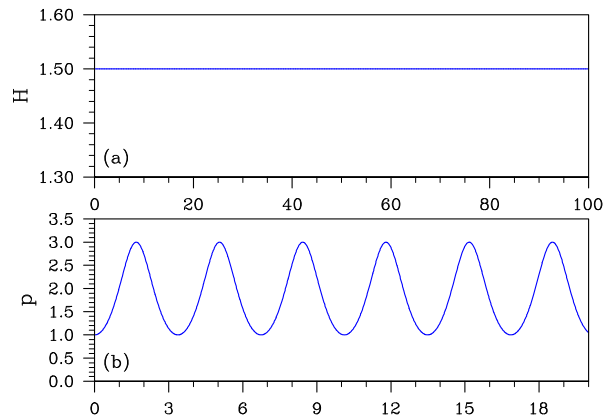

**Figure 2.** The direct simulation of Eq. (3) by the 20-order PMT method, with a step-size of $h = 0.01$: (a) Hamiltonian *H*, (b) variable *p* versus time (the first 20 time units).

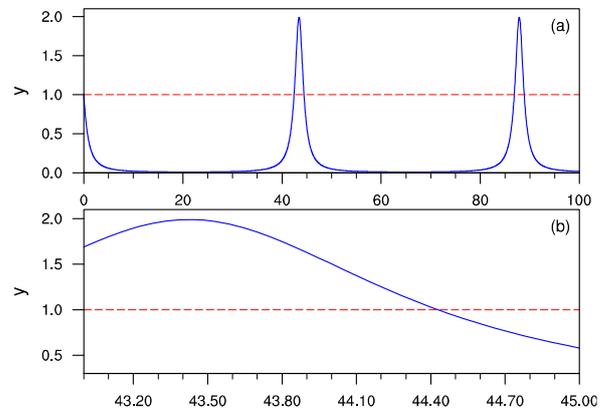

**Figure 3.** Demonstration of the FPA method to obtain the period for Eq. (2): (a) the numerical result of variable *y* by the PMT method in the interval [0,100]; and (b) variable *y* in the enlarged interval of [43,45].